# Comment: The Place of Death in the Quality of Life

**Paul R. Rosenbaum**

## 1. MISTAKING AN OUTCOME FOR A COVARIATE

Donald Rubin's lucid discussion of censoring by death comments on several issues: he warns against mistakes, describes obstacles to inference that might be surmounted within a given investigation, and discusses barriers to inference that direct attention to new data from outside the current investigation. Censoring by death creates outcomes that are defined only contingently, such as quality of life defined only for survivors. If the contingency is an outcome of treatment—if survival could be affected by the treatment—then, as Rubin demonstrates, it is a serious analytical mistake to act as if the contingency were a covariate, a variable unaffected by treatment, when studying the effect of the treatment on the contingently defined outcome. This is one instance of a family of interlinked errors in which an analysis uses an outcome of treatment as if it were a covariate measured before treatment. Other instances in this same family are adjusting for an outcome as if it were a covariate (Rosenbaum, 1984), or attempting to define an interaction effect between a treatment and an outcome of treatment (Rosenbaum, 2004). One of the several advantages of defining outcomes of treatment as comparisons of potential responses under alternative treatments (Neyman, 1923; Rubin, 1974) is that it becomes difficult to make these mistakes: outcomes exist in several versions depending upon the treatment, whereas covariates exist in a single version.

Figure 1 depicts the mistake Rubin warns against. It is a simulated randomized experiment, with $N =$

*Paul R. Rosenbaum is Robert G. Putzel Professor of Statistics, The Wharton School, University of Pennsylvania, Philadelphia, Pennsylvania 19104-6340, USA e-mail: rosenbap@wharton.upenn.edu.*

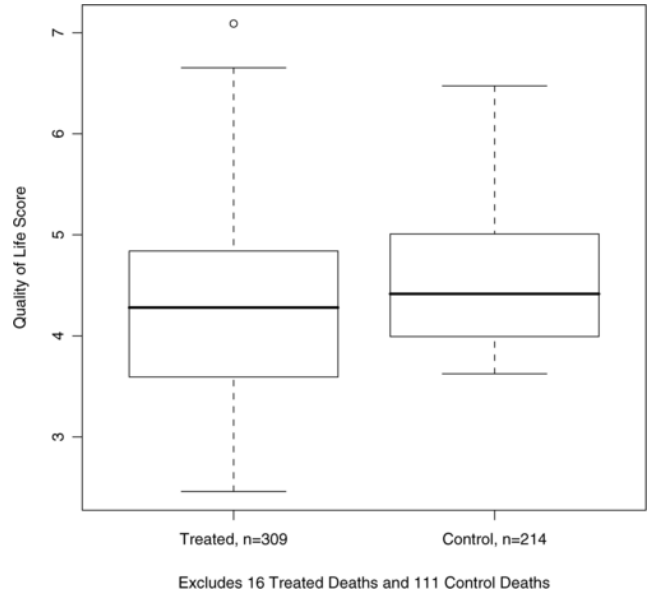

FIG. 1. *Comparison of quality of life.*

650 subjects, of whom $n = 325$ were randomized to treatment where 16 died, and $m = 325$ were randomized to control, where 111 died, and Figure 1 depicts quality of life scores for survivors. Beginning with the structure as Rubin develops it, I will propose a somewhat different analysis. In Section 2 notation describes a completely randomized experiment of the type depicted in Figure 1, with censoring by death but without covariates; then Section 3 proposes a method of analysis that separates empirical evidence of treatment effect from diverging patient preference orderings of death and various qualities of life.

## 2. CENSORING BY DEATH IN A COMPLETELY RANDOMIZED EXPERIMENT

There is a finite population of $N$ subjects, $i = 1, \ldots, N$, who have given informed consent to be randomized to receive either the treatment condition or the control condition, where subject $i$ would exhibit response $r_{Ti}$ under treatment or response $r_{Ci}$ under control. Write $\mathcal{R}$ for $\{(r_{Ti}, r_{Ci}), i = 1, \ldots, N\}$ for the potential responses of the $N$ subjects, which are fixed features of the finite population of $N$ subjects. Of the $N$ subjects, a fixed number, $n$, with $1 \leq n < N$, are picked at random for treatment, denoted $Z_i = 1$, the remaining $m = N - n$ receiving control, denoted $Z_i = 0$, so that $n = \sum_{i=1}^{N} Z_i$, and all $\binom{N}{n}$ possible treatment assignments $\mathbf{Z} = (Z_1, \ldots, Z_N)^T$ have the same probability $\binom{N}{n}^{-1}$. The response, $R_i$, actually observed from $i$ is $r_{Ti}$ if $Z_i = 1$ or $r_{Ci}$ if $Z_i = 0$, and the observed data are $\mathcal{O} = \{(R_i, Z_i), i = 1, \ldots, N\}$. Here $\mathcal{R}$ is fixed but $\mathcal{O}$ is random, and the distribution of $\mathcal{O}$ is created from $\mathcal{R}$ by the known probability distribution used in the random assignment of treatments. The task of inference in a completely randomized experiment is to say something about the effects caused by the treatment, $\mathcal{R}$, from the observed data, $\mathcal{O}$, and the known distribution of treatment assignments. This commonplace description of a randomized experiment is found, for instance, in Welch (1937), and it merged certain ideas from Fisher (1935) about randomization inference and certain ideas from Neyman (1923) about treatment effects.

Following Rubin's approach, I will understand "censoring by death" to mean that the response is a numerical measure of "quality of life" at a particular time, say a year, after treatment, taking values in a subset $\mathcal{Q}$ of the real line, but the measure is not defined if the subject has "died" before that time, in which case the letter "$D$" appears in place of the numerical measure, so $(r_{Ti}, r_{Ci})$ could be a pair of numbers, a $D$ paired with a number, a number paired with a $D$, or a pair of $D$'s. The mistake mentioned in Section 1 consists in setting aside the $D$'s when studying quality of life, and as Rubin's discussion makes very clear, setting aside the deaths means not estimating the effects of the treatment on quality of life.

It is sometimes the case that deaths can be compared ordinally to various qualities of life, even though numerical comparisons are not possible; that is, $\mathcal{Q} \cup \{D\}$ may be a totally ordered set, with strict inequality $\prec$ and with equality-or-inequality $\preceq$, but the elements of $\mathcal{Q} \cup \{D\}$ cannot be manipulated arithmetically to yield averages or expected values. One common view, perhaps the default view, might order death as inferior to any quality of life, and that view might have such diverse sources as religious teachings or the very different observation that a living person can end his or her life, so remaining alive with a given quality of life reveals a preference for that quality of life over death. This common or default view is, no doubt, not universally held, and a particular person might order death or $D$ as preferable to the lowest or worst qualities of life in $\mathcal{Q}$. The analysis that I will describe can accommodate any total ordering of $\mathcal{Q} \cup \{D\}$; it need not place $D$ below all of $\mathcal{Q}$. Faced with diverse preferences among different patients, one can carry out the proposed analysis with several different placements of $D$ in $\mathcal{Q} \cup \{D\}$, in which case the empirical results of a single experiment might speak differently to different patients, and each patient could select the analysis that corresponds to that patient's own evaluation. This is illustrated in Section 3. When a total ordering of $\mathcal{Q} \cup \{D\}$ is possible, to what extent does it facilitate inferences about the effects caused by treatments? More abstractly, what can be said about treatment effects when outcomes take values in a totally ordered set that lacks algebraic operations?

## 3. THE QUALITY OF LIFE AMID DEATH

In the randomized experiment, we observe $n$ of the $N$ potential responses to treatment and we do not observe $m$ of the $N$ potential responses to treatment, and we observe $m$ of the $N$ potential responses to control, but we do not observe $n$ of the $N$ potential responses to control. Let $R_{T(1)} \preceq R_{T(2)} \preceq \cdots \preceq R_{T(n)}$ denote the ordered, observed responses to treatment, including the $D$'s, for the $n$ treated subjects, $Z_i = 1$, and let $\widetilde{R}_{T(1)} \preceq \widetilde{R}_{T(2)} \preceq \cdots \preceq \widetilde{R}_{T(m)}$ denote the unobserved, ordered responses to treatment for the $m$ control subjects, $Z_i = 0$. In Figure 1, there are 16 $D$'s observed in the treated group, so 16 of the $R_{T(i)}$'s are $D$'s, and if deaths are placed below any quality of life by $\preceq$, then the $R_{T(1)} = \cdots = R_{T(16)} = D$. Similarly, let $R_{C(1)} \preceq R_{C(2)} \preceq \cdots \preceq R_{C(m)}$ denote the ordered, observed responses to control for the $m$ control subjects, $Z_i = 0$, and let $\widetilde{R}_{C(1)} \preceq \widetilde{R}_{C(2)} \preceq \cdots \preceq \widetilde{R}_{C(n)}$ denote the ordered, unobserved responses to control for the $n$ treated subjects, $Z_i = 1$. Note that, although $\mathcal{R}$ is fixed, the $R_{T(i)}$, $\widetilde{R}_{T(j)}$, $R_{C(k)}$ and $\widetilde{R}_{C(\ell)}$ are random variables with distributions created from $\mathcal{R}$ by random assignment of treatments; moreover, the $R_{T(i)}$, $\widetilde{R}_{T(j)}$, $R_{C(k)}$ and $\widetilde{R}_{C(\ell)}$ may be numbers in $\mathcal{Q}$ or the letter $D$.

Fix an $i$, $1 \leq i \leq n$, and consider the bivariate random vector $\mathbf{\Upsilon}_{(i)} = \langle R_{T(i)}, \widetilde{R}_{C(i)} \rangle$. Here, $R_{T(i)}$ is the observed $i$th largest response of the $n$ responses of the $n$ subjects randomly assigned to treatment, and $\widetilde{R}_{C(i)}$ is the unobserved $i$th largest response that would have been observed from these same $n$ subjects had they all received the control instead, and either coordinate of $\Upsilon_{(i)}$ may be a $D$. If $n$ were odd and $i = (n+1)/2$, then $\Upsilon_{(i)} = \langle R_{T(i)}, \widetilde{R}_{C(i)}\rangle$ would compare the median of the $n$ observed responses, including deaths, to treatment among $n$ treated subjects to the median of the $n$ unobserved responses, including deaths, that would have been observed among these same $n$ treated subjects had they all received control instead of treatment. Notice carefully that there may be no individual $i$ with $(r_{Ti}, r_{Ci}) = \langle R_{T(i)}, \widetilde{R}_{C(i)}\rangle$, and the quantity $R_{T(i)} - \widetilde{R}_{C(i)}$ is not generally defined because either $R_{T(i)}$ or $\widetilde{R}_{C(i)}$ may equal $D$.

Because $\widetilde{R}_{C(i)}$ is not observed, $\Upsilon_{(i)}$ too is not observed. An exact, randomization-based confidence set for $\Upsilon_{(i)}$ will now be defined. Recall that $\mathcal{C}_{(i)}$ is a $1-\alpha$ confidence set for an unobserved random vector $\Upsilon_{(i)}$ if (i) $\mathcal{C}_{(i)}$ is a function of the observed data, $\mathcal{O}$, and (ii) $1-\alpha \leq \Pr\{\Upsilon_{(i)} \in \mathcal{C}_{(i)}\}$; see Weiss (1955). Proposition 1 rephrases a result due to Fligner and Wolfe (1976, page 83, B; 1979); see Remark 2 following the proposition. The confidence set for $\Upsilon_{(i)}$ is the observed $R_{T(i)}$ and an interval for $\widetilde{R}_{C(i)}$ formed from two of the observed $R_{C(j)}$'s. Notice that the interval may have one or both endpoints as a $D$.

PROPOSITION 1 (Fligner and Wolfe, 1976, 1979). *If $1 \leq a \leq b \leq m$ are two integers such that*

$$(1) \quad 1-\alpha = \sum_{j=a}^{b} \frac{\binom{m+n-i-j}{m-j}\binom{i+j-1}{j}}{\binom{N}{m}},$$

*then $\mathcal{C}_{(i)} = \{\langle R_{T(i)}, w\rangle : w \in [R_{C(a)}, R_{C(b)}]\}$ is a $1-\alpha$ confidence set for $\Upsilon_{(i)}$.*

REMARK 2. Fligner and Wolfe (1976) derive a prediction interval for an order statistic from a future sample starting from i.i.d. sampling of an infinite population, but it is straightforward to derive their combinatorial result, namely their Corollary 4.1 in Fligner and Wolfe (1976), from random assignment of treatments in a finite population, and from this, the coverage of their prediction interval follows. Specifically, (i) start by assuming the $N$ fixed, ordered responses to control are untied, $r_{C(1)} \prec r_{C(2)} \prec \cdots \prec r_{C(N)}$; (ii) then, $\binom{m+n-i-j}{m-j}\binom{i+j-1}{j}$ of the $\binom{N}{m}$ possible random assignments $\mathbf{Z}$ produce Table 1, yielding (1) in agreement with Corollary 4.1 in Fligner and Wolfe (1976); (iii) finally, note with Fligner and Wolfe (1976, page 84) or by other methods that ties among the $r_{Ci}$ make the prediction interval conservative.

At first, adopt the default view, that places death or $D$ below all qualities of life in $\mathcal{Q}$. Then, in Figure 1, there are $n = m = 325$ subjects in each group, and the $R_{T(i)} = D$ for $i = 1, \ldots, 16$, $R_{C(j)} = D$ for $i = 1, \ldots, 111$. With $i = (n+1)/2 = (325+1)/2 = 163$, the median observed response in the treated group is $R_{T(163)} = 4.19$. With $a = 138$, $b = 189$, expression (1) equals 0.951, and $[R_{C(a)}, R_{C(b)}] = [3.81, 4.16]$, so the 95% confidence set for $\Upsilon_{(i)}$ is $\mathcal{C}_{(i)} = \{\langle 4.19, w\rangle : w \in [3.81, 4.16]\}$. This 95% confidence set excludes the possibility that, taking account of the unequal death rates, the median quality of life score would have been higher for the $n = 325$ treated subjects had they all received the control instead, despite the appearance of Figure 1.

Table 2 gives $\mathcal{C}_{(i)}$ for $i = 41$, 82, 163, 244 and 285, for the eighth's, quartiles and median. Notice that for $i = 82$ for the lower quartile, the 95% confidence set contrasts the observed lower quartile in the treated group, $R_{T(82)} = 3.49$, to an interval $[R_{C(61)}, R_{C(106)}] = [D, D]$, so with 95% confidence the lower quartile of the treated group would have been "death" had all $n$ treated subjects received control.

Consider now a hypothetical patient who views qualities of life greater than or equal to 3.5 as better than death, but qualities below 3.5 as inferior to death. What does the same randomized trial say to such a hypothetical patient with these hypothetical preferences? In Figure 1, there are 66 treated patients and no control patients with qualities below 3.5. As a result, with this placement of $D$ in $\mathcal{Q} \cup \{D\}$, the $R_{C(j)}$ are unchanged, but the $R_{T(i)}$ reflect the new order, with $R_{T(i)} \prec D \prec 3.5$ for $i = 1, \ldots, 66$, $R_{T(i)} = D \prec 3.5$ for $i = 67, \ldots, 82$, and $D \prec 3.5 \preceq R_{T(i)}$ for $i = 83, \ldots, n = 325$. Then $\mathcal{C}_{(41)} = \{\langle 3.23, w \rangle : w \in [D, D]\}$ where $3.23 \prec D$ so, with 95% confidence, the lower eighth is worse if all $n$ treated subjects had received control, but $\mathcal{C}_{(82)} = \{\langle D, w \rangle : w \in [D, D]\}$ so the lower quartiles would be the same, and the remaining three intervals in Table 2 are unchanged. With the default order, treatment appeared superior, but with the hypothetical order, control appears better at the lower eighth and worse at the median.

TABLE 1
*Cross-classification of untied potential responses to control, $r_{Ci}$, by treatment assignment, $Z_i$, dividing at $\widetilde{R}_{C(i)}$*

|  | $\prec \widetilde{R}_{C(i)}$ | $\widetilde{R}_{C(i)}$ | $\succ \widetilde{R}_{C(i)}$ | **Total** |
|---|---|---|---|---|
| Treated | $i-1$ | 1 | $n-i$ | $n$ |
| Control | $j$ | 0 | $m-j$ | $m$ |

TABLE 2
*Inference under the default order: 95% confidence set $\mathcal{C}_{(i)}$ for $\Upsilon_{(i)}$*

| **Quantile** | $i$ | $(a,b)$ | $R_{T(i)}$ | $[R_{C(a)}, R_{C(b)}]$ |
|---|---|---|---|---|
| Lower eighth | 41 | (25, 60) | 3.10 | $[D, D]$ |
| Lower quartile | 82 | (61, 106) | 3.49 | $[D, D]$ |
| Median | 163 | (138, 189) | 4.19 | $[3.81, 4.16]$ |
| Upper quartile | 244 | (221, 266) | 4.79 | $[4.43, 4.94]$ |
| Upper eighth | 285 | (267, 302) | 5.20 | $[4.98, 5.58]$ |

Perhaps there is a correct placement of death, $D$, amid the possible qualities of life, $\mathcal{Q}$, or perhaps not. Certain religious teachings would place $D$ below all of $\mathcal{Q}$, but that view is not universal: Seneca (49 A.D., page 92), wrote: "He will live badly who does not know how to die well." The randomized experiment in Section 2 provides no new insight into the proper placement of $D$ in $\mathcal{Q} \cup \{D\}$. However, for each given placement of $D$ in $\mathcal{Q} \cup \{D\}$, the experiment provides information about how a group of $n$ people will fare under treatment and under control.